\documentclass[12pt]{article}
\usepackage{amsmath, amssymb}

\newcommand{\F}{\noindent}

\newcommand{\SP}{\smallskip}
\newcommand{\MP}{\medskip}
\newcommand{\BP}{\bigskip}

\newcommand{\beq}{\begin{eqnarray}}
\newcommand{\ene}{\end{eqnarray}}

\newcommand{\beqs}{\begin{eqnarray*}}
\newcommand{\enes}{\end{eqnarray*}}

\newcommand{\eq}[1]{(\ref{#1})}

\newcommand{\nom}{\nonumber}

\newcommand{\R}{{{\mathbb R}}}

\newcommand{\HH}{{\cal H}}

\newcommand{\DD}{{\cal D}}
\newcommand{\SSS}{{\cal{S}}}

\newcommand{\BB}{{\cal{B}}}

\newcommand{\pmat}{\begin{pmatrix}}
\newcommand{\emat}{\end{pmatrix}}

\newcommand{\hxi}{{\widehat \xi}}
\newcommand{\heta}{{\widehat \eta}}

\Large

\begin{document}

\newtheorem{df}{Definition}[section]

\newtheorem{thm}{Theorem}[section]

\newtheorem{ass}{Assumption}

\newtheorem{lem}{Lemma}

\newtheorem{pro}{Proposition}[section]

\newtheorem{axm}{Axiom}

\numberwithin{equation}{section}

\rightline{KIMS-2006-06-21}

\BP



\begin{center}
\Large
{\bf Fundamental solution global in time for\\
a class of Schr\"odinger equations with\\
time-dependent potentials}
\vskip18pt

\normalsize
Hitoshi Kitada
\vskip2pt

Graduate School of Mathematical Sciences

University of Tokyo

Komaba, Meguro-ku, Tokyo 153-8914, Japan

e-mail: kitada@ms.u-tokyo.ac.jp
\ 

\vskip10pt

(June 21, 2006)
\end{center}

\BP

\leftskip24pt
\rightskip24pt

\small

\noindent
{\it Abstract}: Fundamental solution for a Schr\"odinger equation with a time-dependent potential of long-range type is constructed. The solution is given as a Fourier integral operator with a symbol uniformly bounded global in time, when measured in natural semi-norms of a symbol class.

\leftskip0pt
\rightskip0pt

\vskip 24pt

\section{Introduction and main result}

\F
We consider a Schr\"odinger operator of the form
\beqs
H(t)=-\frac{1}{2}\Delta+V(t,x)\label{1}
\enes
defined in $\HH=L^2(\R^n)$ $(n\ge 1)$.
Here
$$
\Delta=\sum_{k=1}^n\frac{\partial^2}{\partial x_k^2}
$$
is Laplacian with domain $\DD(H)=H^2(\R^n)$, the Sobolev space of order two, and the time-dependent potential $V(t,x)$ satisfies the following assumption. We use the notation: $\partial_x=(\partial/\partial_{x_1},\cdots,\partial/\partial_{x_n})$, $\partial_x^\alpha=(\partial/\partial_{x_1})^{\alpha_1}\cdots(\partial/\partial_{x_n})^{\alpha_n}$ for a multi-index $\alpha=(\alpha_1,\cdots,\alpha_n)$ with $\alpha_j\ge0$ being an integer, $|\alpha|=\alpha_1+\cdots+\alpha_n$, and $\langle y\rangle=\sqrt{1+|y|^2}$ for $y\in \R^d$ $(d\ge 1)$.
\BP

\F
{\bf Assumption (V)}\ \ $V(t,x)$ is a real-valued $C^\infty$ function of $x\in \R^{n}$ for each $t\in \R$ such that the derivatives $\partial_x^\alpha V(t,x)$ are continuous in $(t,x)\in \R^{n+1}$ for any multi-index $\alpha$ and satisfy the condition:
\SP

\F
There exists a constant $\epsilon$ $(1>\epsilon>0)$ such that for any $\alpha$ with $|\alpha|\ge 1$
\beqs
|\partial_x^\alpha V(t,x)|\le C_\alpha\langle t\rangle^{-|\alpha|-\epsilon}\label{2}
\enes
with some constant $C_\alpha>0$ independent of $x\in\R^n$ and $t\in \R$. 
\BP

Thus $H(t)$ is considered a self-adjoint operator with $\DD(H(t))=H^2(\R^n)$ for each $t\in \R$. 

\pagebreak

Under the assumption we will give a construction of the fundamental solution $U(t,s)$ of the Schr\"odinger equation
\beq
(D_t+H(t))U(t,s)f=0,\quad U(s,s)f=f\label{4}
\ene
in the form of a Fourier integral operator
\beqs
U(t,s)f(x)=\int_{\R^n}\int_{\R^n}e^{i(x\xi-\phi(s,t,y,\xi))}u(t,s,\xi,y)f(y)dy d{\hxi}\label{5}
\enes
for $t\ge s\ge T$ or $t\le s\le -T$ for a sufficiently large $T>0$, where the integral is interpreted as an oscillatory integral (see e.g. \cite{K1}, \cite{K2a}), $D_t=-id/dt$, and $\hxi=(2\pi)^{-n}\xi$. A construction of this form of $U(t,s)$ was given in \cite{K0a} with $\phi(s,t,y,\xi)=(t-s)\xi^2/2+y\xi$. However the amplitude function $u(t,s,\xi,y)$ obtained there is not necessarily uniformly bounded in time $t\in\R$ when measured in semi-norms of the symbol class, although the uniformly boundedness of the family $U(t,s)$ with $t,s \in\R$ in operator norm of $L^2(\R^n)$ has been given. There are other works \cite{F}, \cite{K0aa}, \cite{K2a}, etc. in which fundamental solutions of Schr\"odinger equations are constructed under more general assumptions on the time-dependent potentials or Hamiltonians. However the uniformly boundedness in time of the symbols or the amplitude functions of the fundamental solutions has not been considered. The purpose of the present paper is to give a construction of $U(t,s)$ such that the semi-norms of the symbol or the amplitude function $u(t,s,\xi,y)$ are uniformly bounded globally in time.

\SP

The present type of Schr\"odinger operators have been considered in \cite{K0} and the scattering theory for them has been established. Concrete examples have been listed in Example 1.3 of \cite{K0}. We mention an example which was not stated there. Consider $N$-body Hamiltonian
\beqs
H_N=-\frac{1}{2}\Delta+V(x),\quad\quad V(x)=\sum_{1\le i<j\le N}V_{ij}(x_{ij}),
\enes
where $x_{ij}=x_i-x_j\in \R^d$ $(d\ge 1)$ is a relative coordinate between the $i$-th and $j$-th particles. If the pair potential $V_{ij}(y)$ is a $C^\infty$ long-range potential in $y\in \R^d$ satisfying for any $\alpha$
\beq
|\partial_{y}^\alpha V_{ij}(y)|\le C_\alpha\langle y\rangle^{-|\alpha|-\delta}\label{3}
\ene
with some $\delta$ $(1>\delta>0)$, then the time dependent potential $V(t,x)$ defined below satisfies Assumption (V) with $0<\epsilon<\delta$:
$$
V(t,x)=\sum_{1\le i<j\le N}V_{ij}(x_{ij})\chi(\langle \log\langle t\rangle\rangle x_{ij}/\langle t\rangle).
$$
Here $\chi(x)$ is a $C^\infty$ function of $x\in \R^d$ which satisfies $0\le \chi(x)\le 1$, and $\chi(x)=1$ when $|x|\ge 2$ and $\chi(x)=0$ when $|x|\le 1$.
This potential describes the effective contribution of $V(x)$ when the $N$ particles scatter and go to $N$ pieces as $t\to\pm \infty$. This observation was utilized to show the asymptotic completeness of modified wave operators for two-body Schr\"odinger operators with long-range potentials in \cite{IK}, \cite{K3}, \cite{K4}, etc. The method in the present paper gives a result in $N$-body case. In fact we can prove that the range of the (modified) clustered wave operator $W_{b_N}^{+}$ with $b_N$ being a cluster decomposition into $N$ clusters, is equal to the scattering space $S_{b_N}^1$ under assumption \eq{3}. Here $S_{b_N}^1$ is defined in \cite{K1}, Definition 7.4, as a space of wave functions $f$ with natural asymptotic behavior that $e^{-itH_N}f$ decomposes into $N$ clusters in a linear fashion in time $t$, when $t\to+\infty$. The same holds for the case $t\to-\infty$ with an obvious modification. As every observable particle is a scattering particle (see \cite{K1}, Chapter 4), this result gives the {\it effective asymptotic completeness} of the $N$ body scattering problem for general long-range pair potentials. These will be discussed elsewhere.
\SP

To state our main result, we define $H(t,x,\xi)$ $(x,\xi\in \R^n, t\in\R)$ by
$$
H(t,x,\xi)=\frac{1}{2}|\xi|^2+V(t,x).
$$
We will construct a solution $\phi(s,t,x,\xi)$ $(t\ge s\ge T)$ of Hamilton-Jacobi equation
$$
\partial_s\phi(s,t,x,\xi)+H(s,x,\nabla_x\phi(s,t,x,\xi))=0
$$
with initial condition
$$
\phi(s,s,x,\xi)=x\xi,
$$
when $T>0$ is large.
Let $\SSS$ denote the space of rapidly decreasing functions on $\R^n$, and let $\BB$ be a class of complex-valued $C^\infty$ functions $p(x,\xi,y)$ of $x,\xi,y\in \R^n$ with the semi-norms
$$
|p|_\ell=\max_{|\alpha|+|\beta|+|\gamma|\le \ell}\sup_{x,\xi,y\in\R^n}\bigl|\partial_x^\alpha\partial_\xi^\beta\partial_y^\gamma p(x,\xi,y)\bigr|<\infty
$$
for $\ell=0,1,2,\cdots$.
\MP

Our main result is the following, where $1$ denotes the constant function.

\F
\begin{thm}\label{MainTheorem}\ \ Let Assumption {\rm (V)} be satisfied. Then there exists a constant $T_1>0$ and a symbol $u(t,s,\xi,y)\in \BB$ for $t\ge s\ge T(\ge T_1)$ such that the fundamental solution $U(t,s)$ of the equation {\rm \eq{4}} is written as
\beq
U(t,s)f(x)=\iint e^{i(x\xi-\phi(s,t,y,\xi))}u(t,s,\xi,y)f(y)dyd\hxi\label{6}
\ene
for any initial function $f\in \SSS$. The amplitude function $u(t,s)=u(t,s,\xi,y)$ satisfies for any integer $\ell=0,1,2,\cdots$
\beq
\sup_{t\ge s\ge T}|u(t,s)-1|_\ell\le C_\ell\langle T\rangle^{-\epsilon}\label{7}
\ene
for some constant $C_\ell>0$ independent of $T(\ge T_1>0)$. Thus $U(t,s)$ with $t\ge s\ge T(\ge T_1>0)$ defines a family of uniformly bounded operators on $L^2(\R^n)$. The similar estimate holds for the negative time $t\le s\le -T(\le -T_1<0)$. 
\end{thm}

\section{ Phase function $\phi(s,t,x,\xi)$}

We will construct the phase function $\phi(s,t,x,\xi)$ in this section. For this purpose we consider the Hamilton equation:
\beqs
\left\{
\begin{array}{l}
\frac{dq}{dt}(t,s)=\nabla_\xi H(t,q(t,s),p(t,s))=p(t,s),\\
\frac{dp}{dt}(t,s)=-\nabla_x H(t,q(t,s),p(t,s))=-\nabla_xV(t,q(t,s))
\end{array}
\right.\label{8}
\enes
with initial condition
\beqs
q(s,s)=x,\quad p(s,s)=\xi.\label{9}
\enes

This is a system of ordinary differential equations, and can be solved by successive approximation and we have the following estimates. In the following $I$ denotes the identity operator in the sense appropriate to each context.


\begin{pro}\label{pro1}\ \ There exist constants $T_0>0$ and $C_0>0$ such that the following holds.

{\rm i)} For any $t\ge s\ge T_0$ and $x,\xi\in \R^n$
\beqs
&&\left\{
\begin{array}{l}
|q(s,t,x,\xi)-x|+|q(t,s,x,\xi)-x|\le C_0(t-s)\{\langle s\rangle^{-\epsilon}+|\xi|\},\\
|p(s,t,x,\xi)-\xi|+|p(t,s,x,\xi)-\xi|\le C_0\langle s\rangle^{-\epsilon},
\end{array}
\right.\\
&&\left\{
\begin{array}{l}
|\nabla_x q(s,t,x,\xi)-I|\le C_0\langle s\rangle^{-\epsilon},\\
|\nabla_x q(t,s,x,\xi)-I|\le C_0(t-s)\langle s\rangle^{-1-\epsilon},\\
|\nabla_xp(s,t,x,\xi)|+|\nabla_xp(t,s,x,\xi)|\le C_0\langle s\rangle^{-1-\epsilon},
\end{array}
\right.\\
&&\left\{
\begin{array}{l}
|\nabla_\xi q(s,t,x,\xi)-(s-t)I|\le C_0(t-s)\langle s\rangle^{-\epsilon},\\
|\nabla_\xi p(s,t,x,\xi)-I|\le C_0(t-s)\langle s\rangle^{-1-\epsilon},\\
|\nabla_\xi q(t,s,x,\xi)-(t-s)I|\le C_0(t-s)\langle s\rangle^{-\epsilon},\\
|\nabla_\xi p(t,s,x,\xi)-I|\le C_0\langle s\rangle^{-\epsilon}.
\end{array}
\right.
\enes

{\rm ii)} For any $\alpha, \beta$ with $|\alpha|+|\beta|\ge 2$, there is a constant $C_{\alpha\beta}>0$ such that for any $t\ge s\ge T_0$ and $x,\xi\in \R^n$
\beqs
\left\{
\begin{array}{l}
|\partial_\xi^\alpha\partial_x^\beta q(t,s,x,\xi)|\le C_{\alpha\beta}(t-s)\langle s\rangle^{-\epsilon},\\
|\partial_\xi^\alpha\partial_x^\beta p(t,s,x,\xi)|\le C_{\alpha\beta}\langle s\rangle^{-\epsilon}.
\end{array}
\right.
\enes

{\rm iii)} For any $\alpha$ with $|\alpha|\le 1$, there is a constant $C_\alpha>0$ such that for any $t\ge s\ge T_0$
\beqs
|\partial_\xi^\alpha(q(t,s,x,\xi)-x-(t-s)p(t,s,x,\xi))|\le C_\alpha\min\{\langle t\rangle^{1-\epsilon},(t-s)\langle s\rangle^{-\epsilon}\}.
\enes
\end{pro}

For the details of the proof, see \cite{K0}.

\BP

Take $T_1\ge T_0$ such that $C_0\langle T_1\rangle^{-\epsilon}<1/2$ for the constants $T_0, C_0>0$ in Proposition \ref{pro1}. Then by the proposition, the mapping $T_x(y)=x+y-q(s,t,y,\xi):\R^n\to\R^n$ is a contraction mapping when $t\ge s \ge T\ge T_1$. Thus by the fixed point theorem for contraction mapping, there is a unique $y\in\R^n$ such that $T_x(y)=y$ for each $x\in \R^n$. From this follows that the equation $x=q(s,t,y,\xi)$ has a unique solution $y\in \R^n$ for each $x\in \R^n$, and hence defines a mapping $x\mapsto y=y(s,t,x,\xi)$ from $\R^n$ to $\R^n$. It is shown further by contraction mapping theorem that this is a bijection, and by the above proposition this mapping becomes a $C^\infty$-diffeomorphism. Similarly, we can show that $\eta\mapsto \xi=p(t,s,x,\eta)$ is a $C^\infty$-diffeomorphism and has the inverse $\xi\mapsto \eta(t,s,x,\xi)$, which is also a $C^\infty$-diffeomorphism. Summarizing these and using the estimates in Proposition \ref{pro1}, we have the following.

\begin{pro}\label{pro2} There is a constant $T_1\ge T_0$ such that the mappings $y\mapsto x=q(s,t,y,\xi)$ and $\eta\mapsto \xi=p(t,s,x,\eta)$ have the inverse $C^\infty$-diffeomorphisms $x\mapsto y(s,t,x,\xi)$ and $\xi\mapsto \eta(t,s,x,\xi)$ for $t\ge s\ge T\ge T_1$, and they satisfy the following.

{\rm i)}\ \ $q(s,t,y(s,t,x,\xi),\xi)=x,\quad p(t,s,x,\eta(t,s,x,\xi))=\xi.$

{\rm ii)}\ \  $q(t,s,x,\eta(t,s,x,\xi))=y(s,t,x,\xi),\quad p(s,t,y(s,t,x,\xi),\xi)=\eta(t,s,x,\xi).$

{\rm iii)}\ \ There is a constant $C_1>0$ such that for any $t\ge s\ge T$ and $x,\xi\in \R^n$
\beqs
&&\left\{
\begin{array}{l}
|\eta(t,s,x,\xi)-\xi|\le C_1\langle s\rangle^{-\epsilon},\\
|y(s,t,x,\xi)-x-(t-s)\xi|\le C_1\min\{\langle t\rangle^{1-\epsilon}, (t-s)\langle s\rangle^{-\epsilon}\},
\end{array}
\right.\\
&&\left\{
\begin{array}{l}
|\nabla_x y(s,t,x,\xi)-I|\le C_1\langle s\rangle^{-\epsilon},\\
|\nabla_\xi y(s,t,x,\xi)-(t-s)I|\le C_1\min\{\langle t\rangle^{1-\epsilon},(t-s)\langle s\rangle^{-\epsilon}\},
\end{array}
\right.\\
&&\left\{
\begin{array}{l}
|\nabla_x\eta(t,s,x,\xi)|\le C_1\langle s\rangle^{-1-\epsilon},\\
|\nabla_\xi\eta(t,s,x,\xi)-I|\le C_1\langle s\rangle^{-\epsilon}.
\end{array}
\right.
\enes

{\rm iv)}\ \ For any $\alpha,\beta$ with $|\alpha|+|\beta|\ge 2$, there is a constant $C_{\alpha\beta}>0$ such that for any $t\ge s\ge T$ and $x,\xi\in \R^n$\beqs
\left\{
\begin{array}{l}
|\partial_\xi^\alpha\partial_x^\beta\eta(t,s,x,\xi)|\le C_{\alpha\beta}\langle s\rangle^{-\epsilon},\\
|\partial_\xi^\alpha\partial_x^\beta y(s,t,x,\xi)|\le C_{\alpha\beta}(t-s+1)\langle s\rangle^{-\epsilon}.
\end{array}
\right.
\enes
\end{pro}

We fix $T\ge T_1$ arbitrarily below and define the phase function.

\begin{df}\label{df1}\ \ Let
$$
L(t,x,\xi)=\xi\cdot\nabla_\xi H(t,x,\xi)-H(t,x,\xi)=\frac{1}{2}|\xi|^2-V(t,x)
$$
be the Lagrangian, and set
$$
w(s,t,y,\eta)=y\cdot\eta+\int_t^s L(\tau,q(\tau,t,y,\eta),p(\tau,t,y,\eta))d\tau.
$$
Then for $t\ge s\ge T$, we define
\beqs
\phi(s,t,x,\xi)=w(s,t,y(s,t,x,\xi),\xi).\label{10}
\enes
\end{df}

By a direct calculation we can show that the following holds.

\begin{pro}\label{pro3}\ \ Let $t\ge s\ge T$. Then the above $\phi$ satisfies
\beqs
\begin{array}{l}
\nabla_x\phi(s,t,x,\xi)=\eta(t,s,x,\xi),\\
\nabla_\xi\phi(s,t,x,\xi)=y(s,t,x,\xi).
\end{array}
\enes
Furthermore $\phi$ is a solution of the Hamilton-Jacobi equations
\beqs
\begin{array}{l}
\partial_s\phi(s,t,x,\xi)+H(s,x,\nabla_x\phi(s,t,x,\xi))=0,\\
\partial_t\phi(s,t,x,\xi)-H(t,\nabla_\xi\phi(s,t,x,\xi),\xi)=0,\\
\phi(s,s,x,\xi)=x\xi,
\end{array}
\enes
and the function $\phi$ is uniquely determined by these equations.
\end{pro}

\section{Approximate fundamental solution}

We define for $t\ge s\ge T$ and $f\in\SSS$
\beq
E(t,s)f(x)=\iint e^{i(x\xi-\phi(s,t,y,\xi))}f(y)dyd\hxi.\label{approx}
\ene
Since $f\in\SSS$, this oscillatory integral is justified by using the differential operator
$$
P=\langle \nabla_y\psi\rangle^{-2}(1-i\nabla_y\psi\cdot\nabla_y),\quad \psi(x,\xi,y)=x\xi-\phi(s,t,y,\xi).
$$
Namely using the relation
$$
Pe^{i(x\xi-\phi(s,t,y,\xi))}=e^{i(x\xi-\phi(s,t,y,\xi))},
$$
we integrate by parts in \eq{approx}. From the equality $\nabla_y\phi(s,t,y,\xi)=\eta(t,s,y,\xi)$ and Proposition \ref{pro2}, it follows the estimate
$$
C\langle \xi\rangle\le\langle\nabla_y\psi\rangle\le C'\langle\xi\rangle
$$
for some constant $C, C'>0$. Utilizing this estimate and making integrations by parts we see that the integral in \eq{approx} converges.
It is clear that $E$ satisfies $E(s,s)=I$.

We set for $t\ge s\ge T$ and $f\in\SSS$
\beq
G(t,s)f(x)=-i(D_t+H(t))E(t,s)f(x).\label{G}
\ene

Then we have the following theorem.

\begin{thm}\label{Gestimate}\ \ Let Assumption {\rm (V)} be satisfied.  Then $G(t,s)$ can be written as
\beqs
G(t,s)f(x)=\iint e^{i(x\xi-\phi(s,t,y,\xi))}g(t,s,\xi,y)f(y)dyd\hxi
\enes
for $t\ge s\ge T$ and $f\in \SSS$.
The amplitude $g$ is given by
\beqs
&&g(t,s,\xi,y)=\iint e^{-iy\eta}\sum_{\ell,k=1}^n\int_0^1(\partial_{x_k}\partial_{x_\ell}V)(t,\theta y+\widetilde{\nabla}_\xi\phi(s,t,\xi,y,\xi-\eta))d\theta\\
&&
\quad\quad\quad\quad\quad\quad\quad\quad\quad\quad\quad\quad\times\int_0^1 r(\partial_{\xi_k}\partial_{\xi_\ell}\phi)(s,t,y,\xi-r\eta)dr\ dyd{\heta},
\enes
where
\beqs
\widetilde{\nabla}_\xi\phi(s,t,\xi,y,\eta)=\int_0^1\nabla_\xi\phi(s,t,y,\eta+\theta(\xi-\eta))d\theta.
\enes
Thus by Assumption {\rm (V)}, Propositions {\rm \ref{pro2}} and {\rm \ref{pro3}}, the following estimate holds:
For any $\alpha,\beta$, there is a constant $C_{\alpha\beta}>0$ such that for $t\ge s\ge T$ and $\xi,y\in \R^n$
\beqs
|\partial_\xi^\alpha\partial_y^\beta g(t,s,\xi,y)|\le C_{\alpha\beta}\langle t\rangle^{-1-\epsilon}.
\enes
In particular we have for $t\ge s\ge T$
\beqs
|g(t,s)|_\ell\le C_\ell\langle t\rangle^{-1-\epsilon}
\enes
for some constant $C_\ell$ and all $\ell=0,1,2,\cdots$. Further
\beqs
\Vert E(t,s)\Vert\le C
\enes
and
\beqs
\Vert G(t,s)\Vert\le C\langle t\rangle^{-1-\epsilon}
\enes
for some constant $C>0$ independent of $t\ge s\ge T$. Here $\Vert S\Vert$ denotes the operator norm of a linear operator $S:L^2(\R^n)\to L^2(\R^n)$.
\end{thm}

The proof is done by usual calculus of Fourier integral operators. See \cite{K0} for the details. $E(t,s)$ is called an approximate fundamental solution.

\section{Proof of Theorem \ref{MainTheorem}}

By the result of \cite{K0a} we know that $U(t,s)$ exists under our assumption (V). Our purpose is to show that this $U(t,s)$ is expressed in the form of \eq{6} for $t\ge s\ge T$ when $T>0$ is large and the symbol $u(t,s)$ satisfies \eq{7}. We compute using $(D_t+H(t))U(t,s)=0$ and \eq{G}
\beq
E(t,s)^*U(t,s)&=&I+\int_s^t\frac{d}{d\tau}(E(\tau,s)^*U(\tau,s))d\tau\nom\\
&=&I+\int_s^tG(\tau,s)^*U(\tau,s)d\tau.\label{kinji}
\ene
Here we note that for $f\in \SSS$
\beq
E(t,s)^*E(t,s)f(x)=\iint e^{i(\phi(s,t,x,\xi)-\phi(s,t,y,\xi))}f(y)dyd\hxi.\label{EE}
\ene
We have
\beqs
\phi(s,t,x,\xi)-\phi(s,t,y,\xi)=(x-y){\widetilde\nabla}_y\phi(s,t,x,\xi,y),
\enes
where
\beqs
{\widetilde\nabla}_y\phi(s,t,x,\xi,y)=\int_0^1\nabla_y\phi(s,t,y+\theta(x-y),\xi)d\theta
\enes
is close to $\xi$ by Propositions \ref{pro2} and \ref{pro3}. Then we can make a change of variable from $\xi$ to $\eta$ in \eq{EE} by
\beq
\eta={\widetilde\nabla}_y\phi(s,t,x,\xi,y),\label{nablaphi}
\ene
and get
\beqs
E(t,s)^*E(t,s)f(x)=\iint e^{i(x-y)\eta}J(s,t,x,\eta,y)f(y)dyd{\heta}.
\enes
Here
$$
J(s,t,x,\eta,y)=\bigl|\det \nabla_\eta{\widetilde\nabla_y}\phi^{-1}(s,t,x,\eta,y)\bigr|
$$
is the Jacobian of the inverse mapping $\xi={\widetilde\nabla_y}\phi^{-1}(s,t,x,\eta,y)$ of \eq{nablaphi}. By Propositions \ref{pro2} and \ref{pro3}, $J(s,t,x,\eta,y)$ is close to 1, so that $E(t,s)^*E(t,s)$ is close to the identity operator. We set
\beqs
p(t,s,\xi,y)=\iint e^{-iz\eta}\{1-J(s,t,y+z,\xi-\eta,y)\}dzd\heta,
\enes
and define a pseudodifferential operator $P(t,s)$ by
$$
P(t,s)f(x)=\iint e^{i(x-y)}p(t,s,\xi,y)f(y)dyd\hxi.
$$
Then by utilizing Fourier's inversion formula (see \cite{K2}, Proposition 2.1, or \cite{K2a}, Theorem 1.4), we have 
$$
E(t,s)^*E(t,s)f(x)=f(x)-P(t,s)f(x).
$$
By Propositions \ref{pro2} and \ref{pro3}, and integrations by parts, we get for the symbol $p(t,s)=p(t,s,\xi,y)$
\beqs
|p(t,s)|_\ell\le c_\ell\langle T\rangle^{-\epsilon}\quad(\ell=0,1,2,\cdots)\label{estimate}
\enes
for some constant $c_\ell>0$ independent of $t\ge s\ge T(\ge T_1)$. Let $p_\nu(t,s)$ $(\nu=1,2,3,\cdots)$ be the symbol of the multi-product 
$$
P(t,s)^\nu=\overbrace{P(t,s)\cdots P(t,s)}^{\mbox{\scriptsize{$\nu$ factors}}}
$$
of the pseudodifferential operator $P(t,s)$. Then using the estimate for the symbol of the multi-product of pseudodifferential operators (see \cite{K2}, Theorem 2.2, \cite{K1}, section 5.3, or \cite{K2a}, Theorem 1.8), we have
$$
|p_\nu(t,s)|_\ell\le {\widetilde C}_0^\nu C_\ell\nu^\ell(|p(t,s)|_{\ell+2n_0})^\nu
\le  C_\ell\nu^\ell({\widetilde C}_0 c_{\ell+2n_0}\langle T\rangle^{-\epsilon})^\nu
$$
for some even integer $n_0>n$ and some constants ${\widetilde C}_0, C_\ell>0$.
We remark that the factor $\nu^\ell$ comes from the sum
$$
\sum_{\ell_1+\cdots+\ell_\nu\le \ell} 1=\sum_{j=0}^\ell \binom{\nu+j-1}{j}\le C_\ell\nu^\ell
$$
in the estimation formula of the symbol of multi-product (see \cite{K2}, \cite{K1}, \cite{K2a}).
Thus when $T>0$ is large, the series of symbols
$$
q(t,s)=1+p_1(t,s)+p_2(t,s)+\cdots
$$
converges in the symbol space $\BB$ such that 
$$
\sup_{t\ge s\ge T}|q(t,s)-1|_\ell\le \sum_{\nu=1}^\infty\sup_{t\ge s\ge T}|p_\nu(t,s)|_\ell\le Q_\ell(T,\epsilon,n_0)=2C_\ell{\widetilde C}_0 c_{\ell+2n_0}\langle T\rangle^{-\epsilon}<\infty,
$$
and gives the symbol of the inverse operator
$$
\left(I-P(t,s)\right)^{-1}.
$$
Therefore
$$
E(t,s)(I-P(t,s))^{-1}
$$
gives the inverse operator $(E(t,s)^*)^{-1}$ of $E(t,s)^*$. Using the product formula of a Fourier integral operator and a pseudodifferential operator (\cite{K2}, Theorem 3.1, the adjoint of the equation (3.1)-a) there), we have from these the expression
$$
(E(t,s)^*)^{-1}f(x)=\iint e^{i(x\xi-\phi(s,t,y,\xi))}{\tilde e}(t,s,\xi,y)f(y)dyd\hxi,
$$
where
$$
{\tilde e}(t,s,\xi,y)=\iint e^{-i(y-z)\eta}q(t,s,\eta+{\widetilde\nabla}_y\phi(s,t,y,\xi,z),y)dzd\heta
$$
satisfies by integrations by parts and Propositions 2 and 3
\beq
\sup_{t\ge s\ge T}|{\tilde e}(t,s)-1|_\ell\le \sup_{t\ge s\ge T}|q(t,s)-1|_{\ell+2n_0}\le Q_{\ell+2n_0}(T,\epsilon,n_0) <\infty\label{tildee}
\ene
for all $\ell=0,1,2,\cdots$. (See \cite{K2}. The argument above is an easy modification of the one of \cite{K2}, Theorem 3.4.)

Multiplying \eq{kinji} by $(E(t,s)^*)^{-1}$ on the left, we get
\beqs
U(t,s)=(E(t,s)^*)^{-1}\left(I+\int_s^tG(\tau,s)^*U(\tau,s)d\tau\right).
\enes
Applying this expression to $U(\tau, s)$ on the right hand side, and iterating the process, we get with writing $t=\tau_0$
\beq
U(t,s)&=&(E(t,s)^*)^{-1}\Biggl(I+\sum_{\nu=1}^\infty\int_s^{\tau_0}\int_s^{\tau_1}\cdots\int_s^{\tau_{\nu-1}}\nom\\
&&\quad\quad \times G(\tau_1,s)^*(E(\tau_1,s)^*)^{-1}\cdots G(\tau_\nu,s)^*(E(\tau_\nu,s)^*)^{-1}d\tau_\nu\cdots d\tau_1\Biggr).\label{U}
\ene
We can show that $R(\tau,s)=G(\tau,s)^*(E(\tau,s)^*)^{-1}$ $(\tau\ge s\ge T)$ is written as a pseudodifferential operator in a way similar to the above for $E(t,s)^*E(t,s)$. Namely we have
$$
R(\tau,s)f(x)=\iint e^{i(x-y)\eta}{\tilde r}(\tau,s,x,\eta,y)f(y)dyd{\heta},
$$
where
$$
{\tilde r}(\tau,s,x,\eta,y)=g(\tau,s,x,{\widetilde\nabla_y}\phi^{-1}(s,\tau,x,\eta,y))
{\tilde e}(\tau,s,{\widetilde\nabla_y}\phi^{-1}(s,\tau,x,\eta,y),y)J(s,\tau,x,\eta,y).
$$
By \cite{K2}, Proposition 2.1, or \cite{K2a}, Theorem 1.4, we can further rewrite
$$
R(\tau,s)f(x)=\iint e^{i(x-y)\xi}r(\tau,s,\xi,y)f(y)dyd\hxi,
$$
where
$$
r(\tau,s,\xi,y)=\iint e^{-iz\eta}{\tilde r}(\tau,s,y+z,\xi-\eta,y)dzd{\heta}.
$$
Integrations by parts with respect to the variables $z$ and $\eta$ on the right hand side show that $r(\tau,s,\xi,y)$ satisfies for some constants $C_\ell,C'_\ell>0$
\beqs
|r(\tau,s)|_\ell\le C_\ell|{\tilde r}(\tau,s)|_{\ell+2n_0}\le C'_\ell|g(\tau,s)|_{\ell+2n_0}|{\tilde e}(\tau,s)|_{\ell+2n_0}.
\enes
This with Theorem \ref{Gestimate} and \eq{tildee} yields
\beqs
|r(\tau,s)|_\ell\le b_\ell\langle \tau\rangle^{-1-\epsilon}
\enes
for some constant $b_\ell>0$.
Using the estimate for the symbol of the multi-product of pseudodifferential operators again, we see that the symbol $k(t,s)=k(t,s,\xi,y)$ of the pseudodifferential operator
$$
K(t,s)=K(t,s,D_x,X')=I+\sum_{\nu=1}^\infty\int_s^{\tau_0}\int_s^{\tau_1}\cdots\int_s^{\tau_{\nu-1}}
R(\tau_1,s)\cdots R(\tau_\nu,s)~d\tau_\nu\cdots d\tau_1
$$
satisfies the estimate
\beq
|k(t,s)-1|_\ell&\hskip-4pt\le&\hskip-4pt \sum_{\nu=1}^\infty {\widetilde C}_0^\nu C_\ell\nu^\ell\int_s^{\tau_0}\int_s^{\tau_1}\cdots\int_s^{\tau_{\nu-1}}b_{\ell+2n_0}\langle\tau_1\rangle^{-1-\epsilon}\cdots
b_{\ell+2n_0}\langle\tau_{\nu}\rangle^{-1-\epsilon}d\tau_\nu\cdots d\tau_1\nom\\
&\hskip-4pt\le&\hskip-4pt C_\ell\sum_{\nu=1}^\infty \nu^\ell \left({\widetilde C}_0b_{\ell+2n_0}\epsilon^{-1}\langle s\rangle^{-\epsilon}\right)^\nu\label{sum}
\ene
for some constants ${\widetilde C}_0, C_\ell >0$. The right hand side of \eq{sum} converges and is bounded by a finite constant $K_\ell(T,\epsilon,n_0)=2 C_\ell{\widetilde C}_0b_{\ell+2n_0}\epsilon^{-1}\langle T\rangle^{-\epsilon} >0$ independent of $t\ge s(\ge T)$ when $T>0$ is large.

Using the product formula of a Fourier integral operator and a pseudodifferential operator again, we see from these and \eq{U} that $U(t,s)=(E(t,s)^*)^{-1}K(t,s)$ has the form \eq{6} of a Fourier integral operator and the symbol $u(t,s,\xi,y)$ of $U(t,s)$ is given by
$$
u(t,s,\xi,y)=\iint e^{-i(y-z)\eta}{\tilde e}(t,s,\xi,z)k(t,s,\eta+{\widetilde\nabla}_y\phi(s,t,y,\xi,z),y)dzd\heta.
$$
Propositions 2 and 3, integrations by parts, and \eq{tildee} now yield the estimate
\beqs
|u(t,s)-1|_\ell&\le& C'_\ell \left\{\left(1+|{\tilde e}(t,s)-1|_{\ell+2n_0}\right)\left(1+|k(t,s)-1|_{\ell+2n_0}\right)-1\right\}\\
&\le& C'_\ell\left\{\left(1+Q_{\ell+4n_0}(T,\epsilon,n_0)\right)\left(1+ K_{\ell+2n_0}(T,\epsilon,n_0)\right)-1\right\}
\enes
for some constant $C'_\ell>0$ independent of $t\ge s\ge T(\ge T_1>0)$. The proof of Theorem \ref{MainTheorem} is complete.

\end{document}